\documentclass[12pt,a4paper]{article}
\usepackage{indentfirst}
\usepackage{amsmath,amssymb,amsfonts,amsthm,graphics}
\usepackage{makeidx}
\usepackage{amsmath,amssymb,amsfonts,amsthm,graphics,graphicx}
\usepackage{makeidx}
\usepackage{color}
\usepackage{ifpdf}
\usepackage{subfigure}

\newtheorem{thm}{Theorem}

\newtheorem{prob}{Problem}
\newtheorem{lem}{Lemma}
\newtheorem{pro}{Proposition}
\newtheorem{cor}{Corollary}
\newtheorem{fact}{Fact}

\theoremstyle{definition}
\newtheorem{defn}{Definition}

\def\-{\mbox{--}}

\newtheorem{remark}{Remark}

\def\pf{\noindent {\it Proof.} }
\begin{document}

\title{\Large\bf Total proper connection of graphs\footnote{Supported by NSFC No.11371205 and 11531011, and PCSIRT.} }
\author{\small Hui Jiang, Xueliang Li, Yingying Zhang\\
\small Center for Combinatorics and LPMC-TJKLC\\
\small Nankai University, Tianjin 300071, China\\
\small E-mail: jhuink@163.com; lxl@nankai.edu.cn;\\
\small zyydlwyx@163.com}
\date{}
\maketitle
\begin{abstract}

A graph is said to be {\it total-colored} if all the edges and the vertices of the graph is colored. A path in a
total-colored graph is a {\it total proper path} if $(i)$ any two adjacent edges on the path differ in color,
$(ii)$ any two internal adjacent vertices on the path differ in color, and $(iii)$ any internal vertex of the path
differs in color from its incident edges on the path. A total-colored graph is called {\it total-proper connected}
if any two vertices of the graph are connected by a total proper path of the graph. For a connected graph $G$,
the {\it total proper connection number} of $G$, denoted by $tpc(G)$, is defined as the smallest number of colors
required to make $G$ total-proper connected. These concepts are inspired by the concepts of proper connection
number $pc(G)$, proper vertex connection number $pvc(G)$ and total rainbow connection number $trc(G)$ of a
connected graph $G$. In this paper, we first determine the value of the total proper connection number $tpc(G)$
for some special graphs $G$. Secondly, we obtain that $tpc(G)\leq 4$ for any $2$-connected graph $G$ and give examples
to show that the upper bound $4$ is sharp. For general graphs, we also obtain an upper bound for $tpc(G)$. Furthermore,
we prove that $tpc(G)\leq \frac{3n}{\delta+1}+1$ for a connected graph $G$ with order $n$ and minimum degree $\delta$.
Finally, we compare $tpc(G)$ with $pvc(G)$ and $pc(G)$, respectively, and obtain that $tpc(G)>pvc(G)$ for any nontrivial
connected graph $G$, and that $tpc(G)$ and $pc(G)$ can differ by $t$ for $0\leq t\leq 2$.

{\flushleft\bf Keywords}: total-colored graph, total proper connection, dominating set

{\flushleft\bf AMS subject classification 2010}: 05C15, 05C40, 05C69, 05C75.
\end{abstract}

\section{Introduction}

In this paper, all graphs considered are simple, finite and undirected. We refer to the book \cite{B} for
undefined notation and terminology in graph theory. A path in an edge-colored graph is a {\it proper path}
if any two adjacent edges differ in color. An edge-colored graph is {\it proper connected} if any two
vertices of the graph are connected by a proper path of the graph. For a connected graph $G$, the
{\it proper connection number} of $G$, denoted by $pc(G)$, is defined as the smallest number of colors
required to make $G$ proper connected. Note that $pc(G)=1$ if and only if $G$ is a complete graph.
The concept of $pc(G)$ was first introduced by Borozan et al. \cite{BFG} and has been well-studied recently.
We refer the reader to \cite{ALL,GLQ,LLZ,LWY} for more details.

As a natural counterpart of the concept of proper connection, the concept of proper vertex connection
was introduced by the authors \cite{JLZZ}. A path in a vertex-colored graph is a {\it vertex-proper path}
if any two internal adjacent vertices on the path differ in color. A vertex-colored graph is
{\it proper vertex connected} if any two vertices of the graph are connected by a vertex-proper
path of the graph. For a connected graph $G$, the {\it proper vertex connection number} of $G$, denoted
by $pvc(G)$, is defined as the smallest number of colors required to make $G$ proper vertex connected.
Especially, set $pvc(G)=0$ for a complete graph $G$. Moreover, we have $pvc(G)\geq 1$ if $G$ is a
noncomplete graph.

Actually, the concepts of the proper connection and proper vertex connection were motivated from the concepts
of the rainbow connection and rainbow vertex connection. For details about them we refer to a book \cite{LiS} and
a survey paper \cite{LiSS}. Here we only state the concept of the total rainbow connection of graphs, which was
introduced by Liu et al. \cite{LMS} and also studied in \cite{JLZ, S}. A graph is {\it total-colored} if all the
edges and vertices of the graph are colored. A path in a total-colored graph is a {\it total
rainbow path} if all the edges and internal vertices on the path differ in color. A total-colored graph is
{\it total rainbow connected} if any two vertices of the graph are connected by a total rainbow path of the graph.
For a connected graph $G$, the {\it total rainbow connection number} of $G$, denoted by $trc(G)$, is defined as the
smallest number of colors required to make $G$ total rainbow connected. Motivated by the concept of the total rainbow
connection, now for the proper connection and proper vertex connection we introduce the concept of the total proper
connection. A path in a total-colored graph is a {\it total proper path} if $(i)$ any two adjacent edges on the path
differ in color, $(ii)$ any two internal adjacent vertices on the path differ in color, and $(iii)$ any internal vertex of
the path differs in color from its incident edges on the path. A total-colored graph is {\it total proper connected}
if any two vertices of the graph are connected by a total proper path of the graph. For a connected graph $G$,
the {\it total proper connection number} of $G$, denoted by $tpc(G)$, is defined as the smallest number of colors
required to make $G$ total proper connected. It is easy to obtain that $tpc(G)=1$ if and only if $G$ is a complete
graph, and $tpc(G)\geq 3$ if $G$ is not complete. Moreover,
$$tpc(G)\geq \max\{pc(G),pvc(G)\}.\ \ \ \  \ \ \ \ \ \ \ \ \ \  \ \ \ \ \ \ \ \ (*)$$
We can also extend the definition of the total proper connection to that of the total proper $k$-connection $tpc_k(G)$ in a similar
way as the definitions of the proper $k$-connection $pc_k(G)$, proper vertex $k$-connection $pvc_k(G)$ and total rainbow $k$-connection $trc_k(G)$,
which were introduced by Borozan et al. in \cite{BFG}, the present authors in \cite{JLZZ} and Liu et al. in \cite{LMS}, respectively.
However, one can see that when $k$ is larger very little have been known. Almost all known results are on the case for $k=1$.
So, in this paper we only focus our attention on the total proper connection $tpc(G)$ of graphs, i.e., $tpc_k(G)$ for the case $k=1$.

The rest of this paper is organized as follows: In Section $2$, we mainly determine the value of $tpc(G)$ for some special graphs, and
moreover, we present some preliminary results. In Section $3$, we obtain that $tpc(G)\leq 4$ for any $2$-connected graph $G$ and
give examples to show that the upper bound $4$ is sharp. For general graphs, we also obtain an upper bound for $tpc(G)$. In Section $4$,
 we prove that $tpc(G)\leq \frac{3n}{\delta+1}+1$ for a connected graph $G$ with order $n$ and minimum degree $\delta$. In Section $5$,
 we compare $tpc(G)$ with $pvc(G)$ and $pc(G)$, respectively, and obtain that $tpc(G)>pvc(G)$ for any nontrivial
connected graph $G$, and that $tpc(G)$ and $pc(G)$ can differ by $t$ for $0\leq t\leq 2$.

\section{Preliminary results}

In this section, we present some preliminary results on the total proper connection number and determine the value of $tpc(G)$
when $G$ is a nontrivial tree, a complete bipartite graph and a complete multipartite graph.

\begin{pro}\label{pro1} If $G$ is a nontrivial connected graph and $H$ is a connected spanning subgraph of $G$, then $tpc(G)\leq tpc(H)$.
In particular, $tpc(G)\leq tpc(T)$ for every spanning tree $T$ of $G$.
\end{pro}

\begin{pro}\label{pro2} Let $G$ be a connected graph of order $n\geq 3$ that contains a bridge. If $b$ is the maximum number of bridges
incident with a single vertex in $G$, then $tpc(G)\geq b+1$.
\end{pro}

Let $\Delta(G)$ denote the maximum degree of a connected graph $G$. We have the following.

\begin{thm}\label{thm1} If $T$ is a tree of order $n\geq 3$, then $tpc(T)=\Delta(T)+1$.
\end{thm}

\pf Since each edge in $T$ is a bridge, we have $tpc(T)\geq \Delta(T)+1$ by Proposition \ref{pro2}.
Now we just need to show that $tpc(T)\leq\Delta(T)+1$. Let $v$ be the vertex with maximum degree $\Delta(T)$
and $N(v)=\{v_1,v_2,\ldots, v_{\Delta(T)}\}$ denote its neighborhood. Take the vertex $v$ as the root of $T$.
Define a total-coloring $c$ of $T$ with $\Delta(T)+1$ colors in the following way: Let $u$ be a vertex in $T$.
If $u=v$, color $(i)$ $v$ and its incident edges with distinct colors from $A=\{1,2,\ldots,\Delta(T),\Delta(T)+1\}$,
and $(ii)$ $v_i$ with the color from $A\backslash\{c(v),c(vv_i)\}$ for $1\leq i\leq \Delta(T)$. If $u\neq v$,
there exists a father of $u$, say $u'$. Let $N(u)=\{u',u_1,u_2,\ldots,u_{d(u)-1}\}$ denote the neighborhood of $u$.
Color the edges $\{uu_j:1\leq j\leq d(u)-1\} $ with distinct colors from $A\backslash\{c(u),c(uu')\}$, and the
vertex $u_j$ with the color from $A\backslash\{c(u),c(uu_j)\}$ for $1\leq j\leq d(u)-1$.

For any two vertices $x_1$ and $x_2$ in $T$, let $P_i$ be a path from $x_i$ to $v$, where $i\in\{1,2\}$. Next we
shall show that there is a total proper path $P$ between $x_1$ and $x_2$. If $P_1$ and $P_2$ are edge-disjoint,
then $P=x_1P_1vP_2x_2$; otherwise, we walk from $x_1$ along $P_1$ to the earliest common vertex, say $y$, and then
switch to $P_2$ and walk to $x_2$, i.e., $P=x_1P_1yP_2x_2$. Thus, $tpc(T)\leq\Delta(T)+1$, and
therefore, $tpc(T)=\Delta(T)+1$. \qed

The consequence below is immediate from Proposition \ref{pro1} and Theorem \ref{thm1}.

\begin{cor}\label{cor1} For a nontrivial connected graph $G$,
$$tpc(G)\leq \min\{\Delta(T)+1:\ T\ is\ a\ spanning\ tree\ of\ G\}.$$
\end{cor}

A {\it Hamiltonian path} in a graph $G$ is a path containing every
vertex of $G$ and a graph having a Hamiltonian path is a {\it traceable graph}. We get the following result.

\begin{cor}\label{cor2} If $G$ is a traceable graph that is not complete, then $tpc(G)=3$.
\end{cor}

Let $K_{m,n}$ denote a complete bipartite graph, where $1\leq m\leq n$. Clearly,
$tpc(K_{1,1})=1$ and $tpc(K_{1,n})=n+1$ if $n\geq 2$. For $m\geq 2$, we have the result below.

\begin{thm}\label{thm12} For $2\leq m\leq n$, we have $tpc(K_{m,n})=3$.
\end{thm}

\pf Let the bipartition of $K_{m,n}$ be $U$ and $V$, where $U=\{u_1,\ldots,u_m\}$ and $V=\{v_1,\ldots,v_n\}$.
Since $K_{m,n}$ is not complete, it suffices to show that $tpc(K_{m,n})\leq 3$. Now we divide our discussion into two cases.

\textbf{Case 1.} $m=2$.

We first give a total-coloring of $K_{m,n}$ with $3$ colors. Color $(1)$  the vertex $u_1$ and the edge $v_1u_2$ with color $1$,
$(ii)$  the vertex $u_2$ and the edge $u_1v_1$ with color $2$, and $(iii)$ all the other edges and vertices with color $3$.
Then we show that there is a total proper path $P$ between any two vertices $u, v$ of $K_{m,n}$. It is clear that $u$ and $v$ are
total proper connected by an edge if they belong to different parts of the bipartition. Next we consider that $u$ and $v$ are in the
same part of the bipartition. For $u,v\in U$, we have $P=uv_1v$. For $u,v\in V$, if one of them is $v_1$, then $P=uu_1v$; otherwise, $P=uu_1v_1u_2v$.

\textbf{Case 2.} $m\geq3$.

Similarly, we first give a total-coloring of $K_{m,n}$ with $3$ colors. Color the vertices and edges of the cycle $u_1v_1u_2v_2u_3v_3u_1$
starting from $u_1$ in turn with the colors $1,2,3$. For $4\leq i\leq n$ and $4\leq j\leq m$, color $(i)$  $u_3v_i$ with color $1$, $(ii)$
$u_jv_1$ with color $2$, and $(iii)$ all the other edges and vertices with color $3$. Now we show that there is a total proper path $P$ between
any two vertices $u,v$ of $K_{m,n}$. It is clear that $u$ and $v$ are total proper connected by an edge if they belong to different parts of the
bipartition. For $u,v\in U\backslash\{u_2,u_3\}$, we have $P=uv_1u_2v_2u_3v_3v$. For $u,v\in V\backslash\{v_1,v_2\}$, we have $P=uu_1v_1u_2v_2u_3v$.
It can be checked that $u$ and $v$ are total proper connected in all other cases.

Therefore, the proof is complete.\qed

Since any complete multipartite graph has a spanning complete bipartite subgraph, we obtain the following corollary.

\begin{cor}\label{cor3} If $G$ is a complete multipartite graph that is neither a complete graph nor a tree, then $tpc(G)=3$.
\end{cor}

\section{Connectivity}

In this section, we first prove that $tpc(G)\leq 4$ for any $2$-connected graph $G$. Also we show that this upper bound
is sharp by presenting a family of a $2$-connected graphs. Finally, we state an upper bound of $tpc(G)$ for general graphs.

Given a colored path $P=v_1v_2\ldots v_{s-1}v_s$ between any two vertices $v_1$ and $v_s$, we denote by $start_{e}(P)$ the
color of the first edge in the path, i.e., $c(v_1v_2)$, and by $end_{e}(P)$ the last color, i.e., $c(v_{s-1}v_s)$. Moreover,
let $start_{v}(P)$ be the color of the first internal vertex in the path, i.e., $c(v_2)$, and $end_{v}(P)$ be the last color,
i.e., $c(v_{s-1})$. If $P$ is just the edge $v_1v_s$, then $start_{e}(P)=end_{e}(P)=c(v_1v_s),,start_{v}(P)=c(v_s)$, and
$end_{v}(P)=c(v_1)$.

\begin{defn}\label{defn1} Let $c$ be a total-coloring of $G$ that makes $G$
total proper connected. We say that $G$ has the strong property if for any pair of vertices
$u,v\in{V(G)}$, there exist two total proper paths $P_1,P_2$ between them (not necessarily disjoint)
such that $(1)$ $c(u)\neq start_{v}(P_i)$ and $c(v)\neq end_{v}(P_i)$ for $i=1,2$, and $(2)$ both
$\{c(u),start_{e}(P_1),start_{e}(P_2)\}$ and $\{c(v),end_{e}(P_1),end_{e}(P_2)\}$ are $3$-sets.
\end{defn}

Let $G$ be a connected graph and $H$ be a spanning subgraph of $G$. We say that $H$ is a {\it spanning minimally $2$-connected subgraph}
of $G$ if the removal of any edge from $H$ would leave $H$ $1$-connected.

\begin{thm}\label{thm2} Let $G$ be a $2$-connected graph. Then $tpc(G)\leq 4$ and there exists a total-coloring of $G$ with $4$
colors such that $G$ has the strong property.
\end{thm}

\pf Let $G^{'}$ be a spanning minimally $2$-connected subgraph of $G$. We apply induction on the number of ears in an ear-decomposition
of $G^{'}$. The base case is that $G^{'}$ is simply a cycle $C_n=v_1v_2\ldots v_nv_{n+1}(=v_1)$. Obviously, $tpc(C_3)=1$ and $tpc(C_n)=3$
for $n\geq 4$. Next define a total-coloring $c$ of $C_n$ with $4$ colors by
\begin{eqnarray}c(v_iv_{i+1})=
\begin{cases}
1, &if\ i \ is\ odd, 1\leq i\leq 2k-1\ for\ n=2k\ or\ n=2k+1\cr 2, &if\ i\ is\ even, 2\leq i\leq n \ for\ n=2k\ or\ n=2k+1\cr 4, &if\ i=2k+1\ for\ n=2k+1\end{cases}
\end{eqnarray}
and
\begin{eqnarray}c(v_i)=
\begin{cases}
3, &if\ i \ is\ odd, 1\leq i\leq 2k-1\ for\ n=2k\ or\ n=2k+1\cr 4, &if\ i\ is\ even, 2\leq i\leq 2k \ for\ n=2k\ or\ n=2k+1\cr 1, &if\ i=2k+1\ for\ n=2k+1.\end{cases}
\end{eqnarray}
Clearly, the total-coloring $c$ makes $G^{'}$ have the strong property.

In an ear-decomposition of $G^{'}$, let $P$ be the last ear with at least one internal vertex since $G^{'}$ is assumed to be minimally $2$-connected.
And denote by $G_1$ the graph after removal of the internal vertices of $P$. Let $u$ and $v$ be the vertices of $P\cap G_1$ and then $P=uu_1u_2\ldots u_pv$.
By induction hypothesis, there exists a total-coloring of $G_1$ with $4$ colors such that $G_1$ is total proper connected with the strong property.
We give such a total-coloring to $G_1$. Then there exist two total proper paths $P_1$ and $P_2$ from $u$ to $v$ such that $(1)$ $c(u)\neq start_{v}(P_i)$
and $c(v)\neq end_{v}(P_i)$ for $i=1,2$, and $(2)$ both $\{c(u),start_{e}(P_1),start_{e}(P_2)\}$ and $\{c(v),end_{e}(P_1),end_{e}(P_2)\}$ are $3$-sets.
Let $A=\{1,2,3,4\}$. Color the edge $uu_1$ with the color from $A\backslash\{c(u),start_e(P_1),start_e(P_2)\}$, and then total-properly color $P$ from $u$ to $v$
so that $c(u_1)\neq c(u),c(u_p)\neq c(v)$ and $c(u_pv)\neq c(v)$. If $c(u_pv)\notin\{end_{e}(P_1),end_{e}(P_2)\}$, it will become clear that this is the easier
case, and so we consider the case that $c(u_pv)\in\{end_{e}(P_1),end_{e}(P_2)\}$ in the following.

Without loss of generality, suppose that $c(u_pv)=end_{e}(P_2)$. We will show that $G^{'}$ is total proper connected with the strong property under this coloring.
For any two vertices of $G_1$, there exist two total proper paths connecting them with the strong property by induction hypothesis. Since $P\cup P_1$ forms
a total proper connected cycle, any two vertices in this cycle also have the desired paths. Assume that $x\in P\backslash\{u,v\}$ and $y\in G_1\backslash P_1$.
Next we will show that there are two total proper paths from $x$ to $y$ with the strong property.

Since $y,u\in G_1$, there exist two total proper paths $P_{u_1}$ and $P_{u_2}$ starting at $y$ and ending at $u$ with the strong property.
Analogously, there exist two total proper paths $P_{v_1}$ and $P_{v_2}$ starting at $y$ and ending at $v$ with the strong property. Since these paths
have the strong property, suppose that $Q_1=xPuP_{u_1}y$ and $Q_2=xPvP_{v_1}y$ are total proper paths. If $end_{e}(Q_1)\neq end_{e}(Q_2)$, then $Q_1$
and $Q_2$ are the desired pair of paths. Thus, assume that $start_{e}(P_{v_1})=start_{e}(P_{u_1})$.

Then there exists a total proper walk $R_1=xPuP_ivP_{v_2}y$ for some $i\in\{1,2\}$ (suppose $i=1$). If $R_1$ is a path, then $R_1$ and $R_2=Q_2$ are the desired two paths.
Otherwise, let $z$ denote the vertex closest to $y$ on $P_{v_2}$ which is in $P_1\cap P_{v_2}$. Now consider the path $R_1^{'}=xPuP_1zP_{v_2}y$.
If $R_1^{'}$ is a total proper path, then $R_1^{'}$ and $R_2$ are the desired two paths, and so we suppose that $end_{e} (uP_1z)=start_{e}(zP_{v_2}y)$.
Since $P_1$ and $P_{v_2}$ are total proper paths, $c(z)\neq start_v(zP_{v_2}y)$, $c(z)\neq start_v(zP_1v)$ and $end_{e}(vP_1z)\neq end_{e}(uP_1z)$.
Then $end_{e}(vP_1z)\neq start_{e}(zP_{v_2}y)$. Let $S_1=xPvP_1zP_{v_2}y$ and $S_2=Q_1$. Obviously, $S_1$ and $S_2$ are two total proper paths.
Note that $end_{e}(zP_{v_2}y)=start_{e}(P_{v_2})\neq start_{e}(P_{v_1})=start_{e}(P_{u_1})$. Thus, $S_1$ and $S_2$ have the strong property.
Since $tpc(G)\leq tpc(G^{'})$ by Proposition \ref{pro1}, we have $tpc(G)\leq 4$ and there exists a total-coloring of $G$ with 4 colors such that
$G$ has the strong property. This completes the proof of Theorem \ref{thm2}. \qed

\begin{figure}[h,t,b,p]
\begin{center}
\scalebox{0.8}[0.8]{\includegraphics{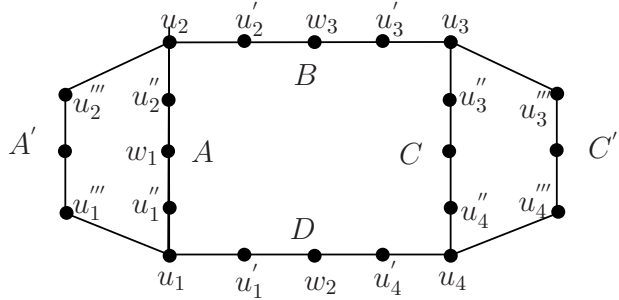}}
\end{center}
\caption{A $2$-connected graph with $tpc(G)=4$.}\label{Fig.1.}
\end{figure}

In order to show that the bound obtained in Theorem \ref{thm2} is sharp, we give a family of $2$-connected graphs $G$ with $tpc(G)=4$ (see Figure \ref{Fig.1.}).

\begin{pro}\label{pro3} Let $G$ be the graph obtained from an even cycle by adding two ears which are as long as their interrupting segments respectively,
such that each segment has $2^{k}$ ($k\geq 2$) edges. Then $tpc(G)=4$.
\end{pro}

Before proving Proposition \ref{pro3}, we give the following fact.

\begin{fact}\label{fact1} Let $C_n=v_1v_2\ldots v_nv_{n+1}(=v_1)$. If there exists a total-coloring of $C_n$ with three colors such that
there are two total proper paths $v_jv_{j+1}\ldots v_{i-1}v_i$ and $v_lv_{l+1}\ldots v_{k-1}v_k$ where $1\leq i<j<k<l\leq n,\mid i-l\mid>1$
and $\mid k-j\mid>1$, then $3\mid n$.
\end{fact}

\noindent {\it Proof of Proposition 3: } Since $tpc(G)\leq 4$ by Theorem \ref{thm2}, we just need to prove that $tpc(G)\neq 3$.
Assume that there is a total-coloring of $G$ with $3$ colors such that $G$ is total proper connected. Label the segments and some
vertices of $G$ as in Figure \ref{Fig.1.}, where $u_i^{'},u_i^{''}$ and $u_i^{'''}$ are the neighbours of the vertex $u_i$ for $i\in\{1,2,3,4\}$.

Firstly, we shall show that the segments $B$ and $D$ are two total proper paths. If one of them is not, say $B$, then there is no total proper path
in $B$ from $u_2$ to $u_3^{'}$ or from $u_3$ to $u_2^{'}$ (say from $u_2$ to $u_3^{'}$). Hence there exists a total proper path through $D$
connecting $u_2$ and $u_3^{'}$, suppose $u_2ADCBu_3^{'}$ (this assumption, as opposed to using any of $A^{'}$ or $C^{'}$, does not lose any generality).
Next we consider the total proper path between $u_1$ and $u_4^{'''}$. Then there must exist a total proper path using the segments $DC^{'}$ or $DCC^{'}$.
If there is a total proper path $u_1DC^{'}u_4^{'''}$, then $c(u_4u_4{''})=c(u_4u_4{'''})$. Thus the total proper path between $u_4^{''}$ and $u_4^{'''}$ is
unique, i.e., $u_4^{''}CC^{'}u_4^{'''}$, and then $c(u_3u_3{'})=c(u_3u_3{'''})$. However, we can not find a total proper path from $u_3^{'}$ to $u_3^{'''}$,
a contradiction. If there is a total proper path $u_1DCC^{'}u_4^{'''}$, then $c(u_3u_3{'})=c(u_3u_3{'''})$. Thus the total proper path connecting $u_3^{'}$
and $u_3^{'''}$ is unique, i.e., $u_3^{'}BCC^{'}u_3^{'''}$. Then $u_3CC^{'}u_3^{'''}$ and $u_4CC^{'}u_4^{'''}$ are two total proper paths in $C\cup C^{'}$
which is an even cycle of length $2^{k+1}$, which contradicts Fact \ref{fact1}. Hence there is no total proper path from $u_1$ to $u_4^{'''}$, a contradiction.
Therefore, the segments $B$ and $D$ are two total proper paths.

Secondly, we will show that at least one of $A$ or $A^{'}$ must be total proper (and similarly, at least one of $C$ or $C^{'}$). Suppose both $A$ and $A^{'}$
are not total proper. Then $u_1$ and $u_2$ are total proper connected by a path through $C$ or $C^{'}$, say $u_1DCBu_2$. However, we can not find a total
proper path connecting $u_1$ and $u_4^{'''}$ in a similar discussion above, which is impossible. Thus, suppose $A$ and $C$ are total proper without loss of generality.

Finally, we know that at least one of the paths $u_1ABCu_4$ and $u_2ADCu_3$ must be not total proper by Fact \ref{fact1}. As we have shown, the only place
which we can not go through is at the intersections, and so assume that the path $w_1ADw_2$  is not total proper, where $w_1\in A\backslash\{u_1,u_1^{''},u_2,u_2^{''}\}$
and $w_2\in D\backslash\{u_1,u_1^{'},u_4,u_4^{'}\}$. In the following, we consider the total proper path $P$ from $w_2$ to $w_1$ and divide our discussion into two cases:

\textbf{Case 1.} $P$ is $w_2DCBAw_1$ or $w_2DCBA^{'}Aw_1$.

Between $w_2$ and $u_4^{'''}$, there must exist a total proper path $P_1$. If $P_1$ is $w_2DC^{'}u_4^{'''}$, then $c(u_4u_4^{''})=c(u_4u_4^{'''})$. Hence,
there is only one total proper path $u_4^{''}CC^{'}u_4^{'''}$ from $u_4^{''}$ to $u_4^{'''}$; otherwise $u_1DCBu_2$ and $u_3BA^{'}Du_4$ are two total proper paths
in the cycle $A^{'}\cup B\cup C\cup D$ for a contradiction. Then it follows that $c(u_3u_3^{'})=c(u_3u_3^{'''})$. Similarly, we can deduce that there is no total proper path
connecting $u_3^{'}$ and $u_3^{'''}$, which is impossible. If $P_1$ is $w_2DCC^{'}u_4^{'''}$, then $c(u_3u_3^{'})=c(u_3u_3^{'''})$. In a similar discussion, we obtain
that the total proper path from $u_3^{'}$ to $u_3^{'''}$ is unique, i.e., $u_3^{'}BCC^{'}u_3^{'''}$. Then $u_3CC^{'}u_3^{'''}$ and $u_4CC^{'}u_4^{'''}$ are two total proper paths
in $C\cup C^{'}$, a contradiction. If $P_1$ is $w_2DA^{'}BC^{'}u_4{'''}$ or $w_2DA^{'}BCC^{'}u_4{'''}$, then $u_1DCBu_2$ and $u_3BA^{'}Du_4$ are two total proper paths
in $A^{'}\cup B\cup C\cup D$, which again contradicts Fact \ref{fact1}.

\textbf{Case 2.} $P$ is $w_2DA^{'}Aw_1$.

Consider the total proper path $P_2$ from $w_2$ to $w_3$, where $w_3\in B\backslash\{u_2,u_2^{'},u_3,u_3^{'}\}$. If $P_2$ is $w_2DCBw_3$, then we can prove that this subcase
could not happen in a similar way as Case 1. If $P_2$ is $w_2DA^{'}Bw_3$, then $c(u_2u_2^{'})=c(u_2u_2^{''})$. From $u_2^{'}$ to $u_2^{''}$, there is only one total proper path
$u_2^{'}BA^{'}Au_2^{''}$ since we can not go through $w_2DCBw_3$. However, $u_2AA^{'}u_2{'''}$ and $u_1A^{'}Aw_1$ are two total proper paths in $A\cup A^{'}$ for a contradiction.

The proof is thus complete. \qed

\begin{remark} Remember that for a $2$-connected graph $G$, we have that the proper connection number $pc(G)\leq 3$; see \cite{BFG}. But, if we consider a $2$-connected bipartite graph $G$, then
we have that $pc(G)=2$. That means that the bipartite property can lower down the number of color by 1. However, from Proposition 3 we see that the bipartite property cannot play a role in general
to lower down the number of colors for the total proper connection number, since the graphs in Proposition 3 are bipartite but their total proper connection numbers reach the upper bound $4$.
\end{remark}

Finally, we prove an upper bound of $tpc(G)$ for general graphs.

\begin{thm}\label{thm3} Let $G$ be a connected graph and $\widetilde{\Delta}$  denote the maximum degree of a vertex which is an endpoint of a bridge in $G$.
Then $tpc(G)\leq\widetilde{\Delta}(G)+1$ if $\widetilde{\Delta}(G)\geq4$ and $tpc(G)\leq4$ otherwise.
\end{thm}

In order to prove Theorem \ref{thm3}, we need a lemma below. Let $R(v)$ denote the set of colors presented on the vertex $v$ and edges incident to $v$.

\begin{lem}\label{lem1} Let $H$ be a graph obtained from a block $B_0$ with $V(B_0)=\{v_1,...,v_n\}$ by adding $t_i(\geq0)$ nontrivial blocks and $s_i(\geq0)$ pendant
edges at $v_i$ for $1\leq i\leq n$. Consider $\widetilde{\Delta}$ as the maximum degree of a vertex which is an endpoint of a bridge in $H$.
Then $tpc(H)\leq \max\{\widetilde{\Delta}(H)+1,4\}$.
\end{lem}

\pf Let $k=\max\{\widetilde{\Delta}(H)+1,4\}$ and $A=\{1,2,...,k\}$. We give a total-coloring $c$ of $H$ using $A$ as follows.

{\bf Step 1.} If $B_0$ is a trivial block, then we give a total-coloring with 3 colors to $B_0$ such that $c(v_1),c(v_2),c(v_1v_2)$ are different from each other;
otherwise, we give a total-coloring with 4 colors to $B_0$ that makes it have the strong property by Theorem \ref{thm2}.
Let $L(v_i)=\{c(v_i),c(v_i)+1,c(v_i)+2,c(v_i)+3\}$ modulo $k$ for $1\leq i\leq n$.

{\bf Step 2.} For $1\leq i\leq n$, if $t_i>0$, then we give a total-coloring with 4 colors from $L(v_i)$ to each uncolored nontrivial block at $v_i$,
denoted by $B^i_j\ (1\leq j\leq t_i)$, that makes each of them have the strong property by Theorem \ref{thm2}; afterwards if $s_i>0$, then color $s_i$ uncolored pendant
edges at $v_i$, denoted by $v_iv^i_1,\ldots,v_iv^i_{s_i}$, with distinct colors from $A\backslash R(v_i)$ and then color each pendant vertex $v^i_m$ using
$A\backslash \{c(v_i),c(v_iv^i_m)\}$ for $1\leq m\leq s_i$.

Next we show that $H$ is total proper connected under the coloring $c$. If $B_0$ is a nontrivial block, then each pair of the vertices in $B_0$ has two total proper paths
between them with the strong property. It will become clear that this is the easier case so we consider the case that $B_0$ is a trivial block. Let $u$ and $w$ be two vertices
of $H$. It is obvious that there exists a total proper path connecting them if both belong to the same block. Suppose that $t_i>0$ and $s_i>0$ for $i=1,2$.
If $u\in\cup^{t_1}_{j=1} B^1_{j}$ and $w\in\cup^{t_2}_{l=1} B^2_l$, then there exist two paths $P_{u_1}$ and $P_{u_2}$ from $u$ to $v_1$ with the strong property.
We know that $uP_{u_j}v_1v_2$ is a total proper path for some $j\in\{1,2\}$ (suppose $j=1$). Similarly, there exists a total proper path $wP_{w_1}v_2v_1$ from $w$ to $v_1$
where $P_{w_1}$ is a total proper path connecting $w$ and $v_2$. Thus, we can find a total proper path $uP_{u_1}v_1v_2P_{w_1}w$ between $u$ and $w$. If $u\in\{v^1_m:1\leq m\leq s_1\}$
and $w\in\{v^2_l:1\leq l\leq s_2\}$, then $uv_1v_2w$ is a total proper path under the coloring $c$. For the other cases, it can be checked that there exists a total proper path
connecting $u$ and $w$ in a similar way. Therefore, $tpc(H)\leq\max\{\widetilde{\Delta}(H)+1,4\}$.  \qed

Now we are ready to prove Theorem \ref{thm3}.

\noindent {\it Proof of Theorem 4:}  Let $B_1,...,B_l$ be the blocks of $G$ and $B(G)$ denote the block graph of $G$ with vertex set $\{B_1,...,B_l\}$. Now,
we consider a breadth-first search tree (BFS-tree) $T$ of $B(G)$ with root $B_1$ and suppose that the blocks have an order $B_1,...,B_l$.
Let $k=\max\{\widetilde{\Delta}(G)+1,4\}$ and $A=\{1,2,...,k\}$. We will give a total-coloring $c$ using $A$ in the following.

We give a total-coloring to $B_1$ and its neighbor blocks of $G$ in a similar way as in Lemma \ref{lem1}. Then we can get that $G$ is total proper connected
if there are no more blocks in $G$. Hence, suppose that there are uncolored blocks in $G$. We extend our coloring from $B_1$ in a Breadth First Search way
until there is no more blocks in $G$, i.e., if $B_i$ has uncolored neighbor blocks, we give a total-coloring to its uncolored neighbor blocks of $G$ in a similar
way as {\bf Step 2}; otherwise, consider $B_{i+1}$.

Now we prove that $G$ is total proper connected. Let $u$ and $w$ be two vertices in $G$. It is obvious that there exists a total proper path between them if
both belong to the same block. Suppose that $u\in B_i$ and $w\in B_j\ (i\neq j)$. Let $P$ denote the path from $B_i$ to $B_j$ in the BFS-tree $T$. Then we can
find a total proper path from $u$ to $w$ traversing the blocks on $P$ under the coloring $c$. Therefore, $tpc(G)\leq\widetilde{\Delta}(G)+1$ if
$\widetilde{\Delta}(G)\geq4$ and $tpc(G)\leq4$ otherwise.
\qed

\section{Minimum degree}

In this section, we prove the following result concerning the minimum degree.

\begin{thm}\label{thm7} Let $G$ be a connected graph of order $n$ with minimum degree $\delta$,
then $tpc(G)\leq \frac{3n}{\delta+1}+1$.
\end{thm}

Given a graph $G$, a set $D\subseteq V(G)$ is called a {\it two-step dominating set} of $G$ if every vertex in $G$ which is not dominated by $D$ has a neighbor
that is dominated by $D$. Moreover, a two-step dominating set $D$ is called a {\it two-way two-step dominating set} if $(a)$ every pendant vertex of $G$ is included in $D$,
and $(b)$ every vertex in $N^2(D)$ has at least two neighbors in $N^1(D)$, where $N^k(D)$ denotes the set of all vertices at distance exactly $k$ from $D$. Further,
if $G[D]$ is connected, $D$ is called a {\it connected two-way two-step dominating set} of $G$.

\begin{lem}\label{lem2}\cite{CDDV} Every connected graph $G$ of order $n\geq 4$ and
minimum degree $\delta$ has a connected two-way two-step
dominating set $D$ of size at most $\frac{3n}{\delta+1}-2$.
\end{lem}

\noindent {\it Proof of Theorem 5:} The proof goes similarly as that of the main result in \cite{LWY} by Li et al.

We are given a connected graph $G$ of order $n$ with minimum degree $\delta$. The assertion can be easily verified for $n\leq3$ and so suppose $n\geq4$.
Let $D$ denote a connected two-way two-step dominating set of $G$ and $k=|D|$. Then we have $k\leq \frac{3n}{\delta+1}-2$ by Lemma \ref{lem2}.
Let $F(x)=\{u:\ u$ is a neighbor of $x$ in $D\}$ for $x\in N^1(D)$ and $F'(y)=\{u: u$ is a neighbor of $y$ in $N^1(D)\}$ for $y\in N^2(D)$.

\textbf{Case 1.} For each vertex $y\in N^2(D)$, its neighbors in $N^1(D)$ has at least one common neighbor in $D$, i.e., $\cap_{x\in F'(y)} F(x)\neq \emptyset$.

\begin{figure}[h,t,b,p]
\begin{center}
\scalebox{0.7}[0.7]{\includegraphics{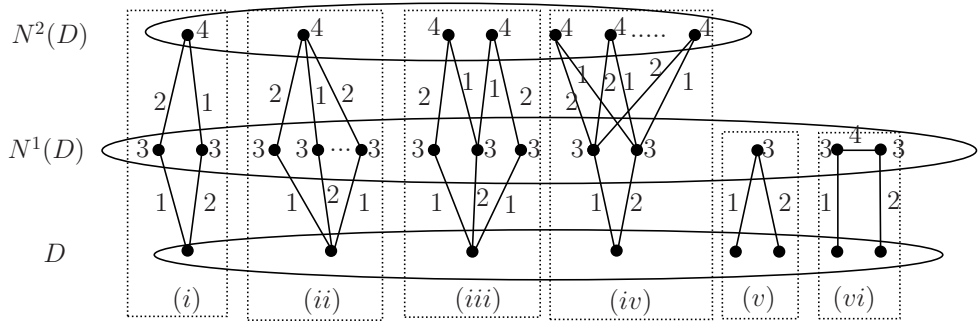}}
\end{center}
\caption{The total-coloring for the spanning subgraph $G_0$ of $G$.}\label{Fig.2.}
\end{figure}

We consider a spanning subgraph $G_0=p(i)\cup q(ii)\cup r(iii)\cup s(iv)\cup t(v)\cup z(vi)\cup G[D]$ of $G$ (see Figure \ref{Fig.2.},
where $p(i)$ denotes the union of $p$ graphs each of which is isomorphic to the graph $(i)$ and similarly for $q(ii),r(iii),s(iv),t(v)$ and $z(vi)$).
Next, we give a total-coloring $c$ to $G_0$ using $\{1,2,...,k,k+1,k+2,k+3\}$. For the edges and vertices of $G[D]$, let $T$ be a spanning tree of $G[D]$.
Then by Theorem 1, $T$ can be total-colored using $\{4,5,...,k,k+1,k+2,k+3\}$ such that for each edge $uv\in E(T)$, the colors of $u,v$ and $uv$ are different from each other.
We color $T$ in such a way and the edges of $G[D]\backslash T$ with any used colors (denote this coloring of $G[D]$ by $c_D$). For the other edges and vertices in $G_0$,
color them as depicted in Figure \ref{Fig.2.}.

Since each pair of vertices $u,w\in D$ has a total proper path $P$ connecting them such that $c(u)\notin\{start_v(P),start_e(P)\}$ and $c(w)\notin\{end_v(P),end_e(P)\}$,
it suffices to show that $G_0$ is total proper connected in the assumption that $\cap_{y\in N^2(D)}\{F(x):\ x\in F'(y)\}=\{w\}$. Take any two vertices $u$ and $v$ in $V(G_0)$.
If $u,v\in N^2(D)$, then $u$ has a neighbor $u'$ in $N^1(D)$ and similarly $v$ has a neighbor $v'$ in $N^1(D)$. Hence, if $c(u'w)\neq c(v'w)$, $uu'wv'v$ is a total proper path;
otherwise, $uu''wv'v$ is a total proper path where $u''$ is another neighbor of $u$ in $N^1(D)$. It is easy to check that $u$ and $v$ are total proper connected in all other cases.

\textbf{Case 2.} There exists one vertex $y\in N^2(D)$ whose neighbors in $N^1(D)$ has no common neighbors in $D$, i.e., $\cap_{x\in F'(y)} F(x)= \emptyset$.

\begin{figure}[h,t,b,p]
\begin{center}
\scalebox{0.7}[0.7]{\includegraphics{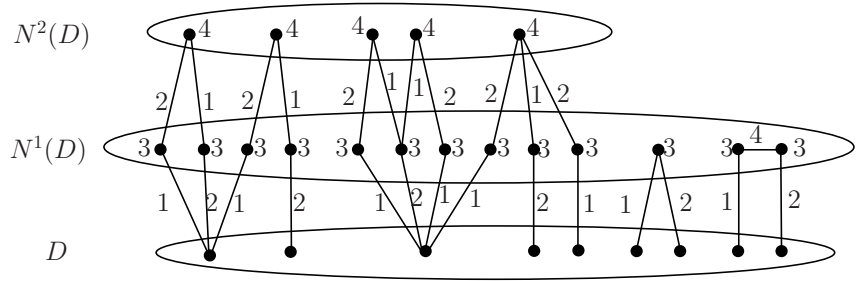}}
\end{center}
\caption{An example for the total-coloring for the spanning subgraph $G_0$ of $G$.}\label{Fig.3.}
\end{figure}

We consider a spanning subgraph $G_0$ of $G$ (see Figure \ref{Fig.3.}). Next, we give a total-coloring $c$ to $G_0$ using $\{1,2,...,k,k+1,k+2,k+3\}$. For the edges and vertices
in $G[D]$, we use the total-coloring $c_D$ as in Case 1. For any vertex $v\in N^2(D)$, color $vx_1$ with color 1 and $x_1u_1$ with color 2 where $x_1\in F'(v)$ and $u_1\in F(x_1)$.
And then color $vx_i$ with color 2 and $x_iu_i$ with color 1 where $x_i\in F'(v)\backslash\{x_1\}$ and $u_i\in F(x_i)$. For any vertex $v\in N^1(D)\backslash \cup_{y\in{N^2}(D)}F^{'}(y)$,
color the edges incident to $v$ as depicted in Figure \ref{Fig.3.}. Moreover, we color the vertices of $N^1(D)$ with color 3 and the vertices of $N^2(D)$ with color 4.

Now we show that $G_0$ is total proper connected. Take any two vertices $u$ and $v$ in $V(G_0)$. If $u,v\in N^2(D)$, there exist two paths $uu'u''$ and $vv'v''$ connecting to $D$,
where $u',v'\in N^1(D)$ and $u'',v''\in D$. Thus, if $u''\neq v''$, $u$ and $v$ are total proper connected by a path $uu'u''Pv''v'v$ where $P$ is a total proper path from $u''$ to
$v''$ in $G[D]$; otherwise, there exists a total proper path connecting $u$ and $v$ in a similar discussion as Case 1. It can be checked that $u$ and $v$ are total proper connected
in all other cases. Therefore, we have $tpc(G)\leq tpc(G_0)\leq\frac{3n}{\delta+1}-2+3=\frac{3n}{\delta+1}+1$ by Proposition \ref{pro1} and Lemma \ref{lem2}. \qed

\section{Compare $tpc(G)$ with $pvc(G)$ and $pc(G)$ }

Let $G$ be a nontrivial connected graph. Recall that $tpc(G)\geq \max\{pc(G), pvc(G)\}$. The question we may ask is, how tight are the inequalities $tpc(G)\geq pc(G)$ and $tpc(G)\geq pvc(G)$ ?
By \cite[Proposition 1 and Theorem 1]{JLZZ}, we have that $(1)$ $pvc(G)=0$ if and only if $G$ is a complete graph, $(2)$ $pvc(G)=1$ if and only if $diam(G)=2$, and $(3)$ $pvc(G)=2$ if and only if $diam(G)\geq 3$. Note that $tpc(G)=1$ if and only if $G$ is a complete graph, and $tpc(G)\geq 3$ if $G$ is not complete. Thus, it follows that $tpc(G)>pvc(G)$.

Next we consider the tightness of the inequality $tpc(G)\geq pc(G)$. Observe that $tpc(G)=pc(G)=1$ if and only if $G$ is a complete graph. Proposition \ref{pro4} below shows that there exists an
example graph $G$ such that $tpc(G)=pc(G)=3$.

\begin{pro}\label{pro4} Let $G$ be the graph obtained from a cycle by adding three ears of length $3$ such that each segment of the cycle has $6t\ (t\geq1)$ edges. Then $tpc(G)=pc(G)=3$.
\end{pro}

\begin{figure}[h,t,b,p]
\begin{center}
\scalebox{0.8}[0.8]{\includegraphics{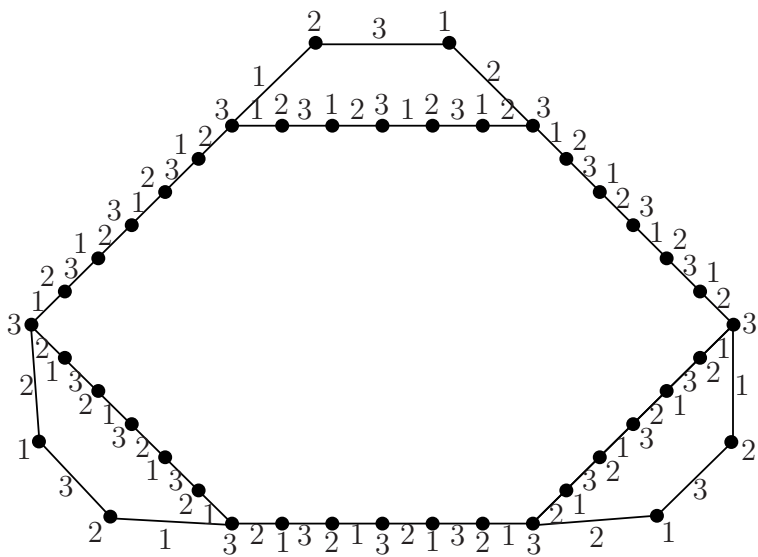}}
\end{center}
\caption{A $3$-coloring of edges and vertices of $G$.}\label{Fig.4.}
\end{figure}

\pf It can be verified that $pc(G)=3$ by \cite[Proposition 3]{BFG}. Thus, it suffices to show that $tpc(G)\leq 3$ by Ineq.$(*)$. A $3$-coloring of edges and vertices of $G$ is shown in Figure \ref{Fig.4.}
to make $G$ total proper connected. Hence, we have that $tpc(G)=pc(G)=3$.
\qed

However, we cannot show whether there exists a graph $G$ such that $tpc(G)=pc(G)=k$ for any $k\geq 4$. Thus, we propose the following problem.

\begin{prob}\label{prob1} For $k\geq4$, does there exist a graph $G$ such that $tpc(G)=pc(G)=k$ ?
\end{prob}

Now we consider the difference between $tpc(G)$ and $pc(G)$. If $T$ is a tree of order $n\geq 3$, then $pc(T)=\Delta(T)$ by \cite[Proposition 2.3]{ALL}
and $tpc(T)=\Delta(T)+1$ by Theorem \ref{thm1}. Hence, $tpc(T)=pc(T)+1$. Moreover, there exists an example graph depicted as in Proposition \ref{pro3}
such that $tpc(G)=pc(G)+2$ since $tpc(G)=4$ and $pc(G)=2$ (we give a $2$-edge-coloring of $G$ by coloring alternately the edges of the segments $A^{'},C^{'}$
and the cycle $A\cup B\cup C\cup D$). However, we have not found any graph $G$ such that $tpc(G)$ and $pc(G)$ can differ by $t\ (t\geq 3)$. Thus, we pose the following problem.

\begin{prob}\label{prob2} For $t\geq3$, does there exist a graph $G$ such that $tpc(G)=pc(G)+t$ ?
\end{prob}

\end{document}